\RequirePackage{ifpdf}
\ifpdf 
\documentclass[pdftex]{sigma}
\else
\documentclass{sigma}
\fi

\newcommand{\be}{\begin{equation}}
\newcommand{\ee}{\end{equation}}
\newcommand{\ba}{\begin{eqnarray}}
\newcommand{\ea}{\end{eqnarray}}
\newcommand{\baa}{\begin{eqnarray*}}
\newcommand{\eaa}{\end{eqnarray*}}
\newcommand{\bb}{\begin{thebibliography}}
\newcommand{\eb}{\end{thebibliography}}

\newcommand{\bi}[1]{\bibitem{#1}}
\newcommand{\lab}[1]{\label{#1}}
\newcommand{\re}[1]{(\ref{#1})}

\numberwithin{equation}{section}

\newcommand{\sn}{\mathop{\rm sn}\nolimits}

\begin{document}

\allowdisplaybreaks

\renewcommand{\PaperNumber}{003}

\FirstPageHeading

\renewcommand{\thefootnote}{$\star$}

\ShortArticleName{Elliptic Biorthogonal Polynomials Connected with
Hermite's Continued Fraction}

\ArticleName{Elliptic Biorthogonal Polynomials Connected\\ with
Hermite's Continued Fraction\footnote{This paper is a contribution
to the Vadim Kuznetsov Memorial Issue ``Integrable Systems and
Related Topics''. The full collection is available at
\href{http://www.emis.de/journals/SIGMA/kuznetsov.html}{http://www.emis.de/journals/SIGMA/kuznetsov.html}}}

\Author{Luc VINET~$^\dag$ and Alexei ZHEDANOV~$^\ddag$}
\AuthorNameForHeading{L.~Vinet and A.~Zhedanov}

\Address{$^\dag$~Universit\'e de Montr\'eal, PO Box 6128, Station
Centre-ville, Montr\'eal QC H3C 3J7, Canada}
\EmailD{\href{mailto:luc.vinet@umontreal.ca}{luc.vinet@umontreal.ca}}

\Address{$^\ddag$~Donetsk Institute for Physics and Technology,
Donetsk 83114, Ukraine}
\EmailD{\href{mailto:zhedanov@yahoo.com}{zhedanov@yahoo.com}}

\ArticleDates{Received October 07, 2006, in f\/inal form December
12, 2006; Published online January 04, 2007}

\Abstract{We study a family of the Laurent biorthogonal
polynomials arising from the Hermite continued fraction for a
ratio of two complete elliptic integrals. Recurrence
coef\/f\/icients, explicit expression and the weight function for
these polynomials are obtained. We construct also a new explicit
example of the Szeg\"{o} polynomials orthogonal on the unit
circle. Relations with associated Legendre polynomials are
considered.}

\Keywords{Laurent biorthogonal polynomials; associated Legendre
polynomials; elliptic integrals}

\Classification{33C45; 42C05}

\begin{flushright}
\it To the memory of Vadim B. Kuznetsov
\end{flushright}

\begin{quote}
\it Vadim Kuznetsov spent a number of years at the Centre de
Recherches Math\'ematiques of the Universit\'e de Montr\'eal in
the mid-90s. This is when we had the privilege to have him as a
colleague and the chance to appreciate his scientif\/ic and human
qualities. His smile, great spirit, penetrating insights and
outstanding scientif\/ic contributions will always be with us.
\end{quote}

\section{Introduction}
We start with well known identity (see e.g. \cite{WW}) \be
\frac{d^2 \sn^n u}{d u^2} = n(n-1) \sn^{n-2} u  - n^2 (1+k^2)
\sn^n u + n(n+1) k^2 \sn^{n+2} u \lab{sn_rec} \ee for the Jacobi
elliptic function $\sn u$ depending on an argument $u$ and a
modulus $k$. This (and similar 11 formulas obtained if one
replaces $sn$ with other elliptic functions) identity goes back to
Jacobi and usually is exploited in order to establish recurrence
relations for elliptic integrals.

Indeed, introduce the following elliptic integrals \be
J_n(x;k)=\int_0^v {k^{2n} \sn^{2n}u \, du}=
    \int_0^x { \frac{k^{2n} t^{2n} dt}{\sqrt{(1-t^2)(1-k^2t^2)}} },
\lab{J_n} \ee where $\sn u =t$ and $x=\sn  v$. Then from
\re{sn_rec} we obtain
\begin{gather} (2n-1) J_n(x;k) -(2n-2)(k^2+1)
J_{n-1}(x;k) + k^2 (2n-3) J_{n-2}(x;k) \nonumber\\
\qquad{}{}= k^{2n-2} x^{2n-3} \sqrt{(1-x^2)(1-k^2x^2)}.
\lab{J_n_rec}
\end{gather}
From this formula we can express $J_{n+1}(x;k)$ in the form
\cite{Her} \be J_{n+1}(x;k) = Q(x;n)  \sqrt{(1-x^2)(1-k^2x^2)} +
A_n J_1(x;k) - B_n J_0(x;k), \lab{JnQ} \ee where $Q(x;n)$ is a
polynomial in $x$ (depending on $n$) and the coef\/f\/icients
$A_n$, $B_n$ satisfy the same recurrence relations as
$J_{n+1}(1;k)$, i.e.
 \begin{gather}
 (2n+1)
A_{n} - 2n(1+k^2) A_{n-1} +(2n-1) k^2 A_{n-2} =0, \nonumber \\
(2n+1) B_{n} - 2n(1+k^2) B_{n-1} +(2n-1) k^2 B_{n-2} =0
\lab{rec_AB}
\end{gather}
 with obvious initial conditions \be A_{-1}=0, \qquad
B_{-1}=-1, \qquad A_0=1, \qquad B_0=0. \lab{ini_AB} \ee Note that
from recurrence relations \re{rec_AB} and initial conditions
\re{ini_AB} it follows that $A_n(k)$ is a~polynomial in $k^2$ of
degree $n$ and $B_n(k)$ is a polynomial in $k^2$ of degree $n$ but
having common factor $k^2$, i.e.\ $B_n=k^2 V_{n-1}(k^2)$, where
$V_{n-1}(z)$ is a polynomial of degree $n-1$ in $z$ for any
$n=1,2,\dots$.

There is an elementary Wronskian-type identity following directly
from \re{rec_AB} \cite{Her}: \be B_n A_{n-1}-A_n
B_{n-1}=\frac{k^{2n}}{2n+1}. \lab{W_AB} \ee

Formula \re{JnQ} allows to reduce calculation of any (incomplete)
elliptic integrals of the form
\[
\int_0^x \frac{P(t) \, dt}{\sqrt{(1-t^2)(1-k^2t^2)}}
\]
(where $P(t)$ is a polynomial) to standard elliptic integrals of
the f\/irst and second kind $J_0(x;k)$ and $J_1(x;k)$. This result
is well known since Jacobi.

In what follows we will denote $J_0(1;k)=K(k)$ and
$J_1(1;k)=J(k)$. Note that $K(k)$ is the standard complete
elliptic integral of the f\/irst kind \cite{WW} and
\[
J(k) = K(k)-E(k)= -k \frac{dE}{dk},
\]
where
\[
E(k)= \int_{0}^1 {   \left(\frac{1-k^2x^2}{1-x^2}\right)^{1/2} \,
dx }
\]
is the complete elliptic integral of the second kind.

We note also an important relation with hypergeometric functions
\cite{WW}:
\begin{gather}
K(k)=\frac{\pi}{2}   {_{2}F_1}(1/2,1/2;1;k^2), \qquad
E(k)=\frac{\pi}{2}\: {_{2}F_1}(-1/2,1/2;1;k^2), \nonumber
\\
J(k)=\frac{\pi k^2}{4}  {_{2}F_1}(3/2,1/2;2;k^2). \lab{KEJ_hyp}
\end{gather}

Hermite in his famous ``Cours d'analyse'' \cite{Her} (see also
\cite{Her1,Her2}) derived a continued fraction connected with a ratio
of two complete elliptic integrals: \be \frac{J(k)}{K(k)} =
{k^2\over\displaystyle 2(1+k^2) - {\strut 9 k^2 \over\displaystyle
4(1+k^2) -{\strut 25 k^2 \over {6(1+k^2) - \cdots {\atop
\displaystyle  }}}}}, \lab{Her_cf} \ee where
\[
K(k)= \int_{0}^1 {   \frac{dx}{\sqrt{(1-x^2)(1-k^2x^2)}}   },
\qquad J(k)= \int_{0}^1 { \frac{k^2x^2 \,
dx}{\sqrt{(1-x^2)(1-k^2x^2)}} }.
\]
The Hermite continued fraction \re{Her_cf} follows directly
from~\re{J_n_rec}. Note that the continued fraction~\re{Her_cf}
belongs to the class of the so-called T-continued fractions (with
respect to the va\-riab\-le~$k^2$) introduced and studied by Thron
(see \cite{JT} for details). Thus perhaps Hermite was the f\/irst
to introduce an explicit example of the $T$-continued fraction. In
what follows we will see that this example gives rise to a class
of polynomials which are biorthogonal on the unit circle.

Rewrite relation \re{JnQ} for $x=1$ in the form \be J(k)/K(k) -
B_n/A_n = \frac{J_{n+1}(k)}{K(k)A_n} \lab{Her_Pade} \ee from which
Hermite concluded that the rational function $B_n/A_n$ is an
approximate expression for the ratio $J(k)/K(k)$ to within terms
of degree $n+1$ in $k^2$.

Hermite also noted that the coef\/f\/icients $A_n$, $B_n$ appeared
as power coef\/f\/icients in the following Taylor expansions:
\begin{gather}
\frac{J_0(x)}{\sqrt{(1-x^2)(1-k^2x^2)}} = A_0 x + A_1 x^3 +
\dots + A_n x^{2n+1} + \cdots, \nonumber \\
\frac{J_1(x)}{\sqrt{(1-x^2)(1-k^2x^2)}} = B_0 x + B_1 x^3 + \dots
+ B_n x^{2n+1} + \cdots. \lab{AB_ser}
\end{gather}
 Thus formulas
\re{AB_ser} can be considered as generating functions for $A_n$,
$B_n$.

\section{Laurent biorthogonal polynomials}

 Introduce the new variable $z=k^2$ and
def\/ine the polynomials $P_n(z)=A_n/\xi_n$ of degree $n$ in $z$,
where
\[
\frac{\xi_n}{\xi_{n+1}}=\frac{2n+3}{2n+2}, \qquad \xi_0=1.
\]
We have \be \xi_n=\frac{n!}{(3/2)_n}, \lab{xi} \ee where
$(a)_n=a(a+1) \cdots(a+n-1)$ is the standard Pochhammer symbol
(shifted factorial).

Then it is seen that $P_n(z) = z^n + O(z^{n-1})$, i.e.\ $P_n(z)$
are monic polynomials. Moreover from~\re{rec_AB} it follows that
$P_n(z)$ satisfy the 3-term recurrence relation \be P_{n+1}(z)
+d_n P_n(z) = z(P_n(z) + b_n P_{n-1}(z)), \lab{rec_LBP}
 \ee where
\be d_n = -1, \quad n=0,1,\dots, \qquad b_n = -
\frac{(n+1/2)^2}{n(n+1)}, \quad n=1,2,\dots \lab{db_Her} \ee with
initial conditions \be P_{-1}(z)=0, \qquad P_0(z)=1. \lab{ini_P}
\ee Def\/ine also the polynomials
\[
P_{n-1}^{(1)}(z) = 2B_n/(z \xi_n ).
\]
 It is then easily verif\/ied that the polynomials $P_n^{(1)}(z)$ are
again $n$-th degree monic polynomials in $z$ satisfying the
recurrence relation \be P_{n+1}^{(1)}(z) + d_{n+1} P_n^{(1)}(z) =
z(P_n^{(1)}(z) + b_{n+1}P_{n-1}^{(1)}(z)). \lab{rec_P1} \ee

The polynomials $P_n(z)$ are Laurent biorthogonal polynomials
(LBP) \cite{HR} and the polynomials~$P_n^{(1)}(z)$ are the
corresponding associated LBP. The recurrence coef\/f\/icients
$b_n$, $d_n$ completely characterize LBP. The nondegeneracy
condition $b_nd_n \ne 0$, $n=1,2,\dots$ \cite{HR,ZheL} obviously
holds in our case.

It is well known that LBP possess the biorthogonality property
\cite{HR}. This means that there exists a family of other LBP
$\hat P_n(z)$ and a linear functional $\sigma$ such that \be
\langle \sigma, P_n(z) \hat P_m(1/z) \rangle = h_m  \delta_{nm},
\lab{bi_LBP} \ee where the normalization constants $h_n$ are
expressed as \cite{HR}
\[
h_n = \prod_{k=1}^n \frac{b_k}{d_k}.
\]

The linear functional $\sigma$ is def\/ined on the space of all
monomials $z^s$ with both positive and negative values of $s$: \be
c_s = \langle \sigma, z^s \rangle, \qquad s=0, \pm 1, \pm 2,
\dots, \lab{moms} \ee where $c_s$ is a sequence of moments (this
sequence is inf\/inite in both directions).

The biorthogonality condition \re{bi_LBP} is equivalent (under the
nondegeneracy condition $h_n \ne 0$) to the orthogonality
relations \be \langle \sigma, P_n(z) z^{-j}\rangle =0, \qquad
j=0,1,\dots,n-1.  \lab{ort_LBP} \ee

Note that in our case from \re{rec_LBP} and \re{db_Her} it follows
that \be P_n(0)=1, \qquad n=0,1,2,\dots. \lab{P01} \ee Introduce
the reciprocal LBP by the formula \cite{VZ1}
 \be  P_n^*(z) = \frac{z^n
P_n(1/z)}{P_n(0)}. \lab{rec_P_def} \ee It appears that $P_n^*(z)$
are again LBP with the recurrence coef\/f\/icients \cite{VZ1} \be
b_n^* = \frac{b_n}{d_n d_{n+1}}, \qquad d_n^* = \frac{1}{d_n}.
\lab{rec_bd} \ee The moments $c_n^*$ of the reciprocal polynomials
are expressed as \cite{VZ1} \be c_n^* = \frac{c_{1-n}}{c_1}.
\lab{mom_rec} \ee

In our case \re{db_Her} we have $b_n^* = b_n$, $d_n^* = d_n.$
Hence the reciprocal polynomials coincide with the initial
ones:{\samepage \be z^n P_n(1/z)=P_n(z). \lab{rec_coin} \ee
Moreover in our case $c_1=d_0=-1$ and thus $c_n^* = -c_{1-n}$.}

The polynomials $\hat P_n(z)$ are the biorthogonal partners with
respect to the polynomials $P_n(z)$. Their moment sequence $\hat
c_n$ is obtained from the initial moment sequence by ref\/lection:
\be \hat c_n = c_{-n}, \qquad n=0, \pm 1, \pm 2, \dots .
\lab{biort_c} \ee Note that under the assumption $c_0=1$ (standard
normalization condition) we have $\hat c_0=1$ as well.

For the biorthogonal partners there is an explicit expression
\cite{HR, VZ1} \be \hat P_n(z) = \frac{z^n P_{n+1}(1/z) - z^{n-1}
P_n(1/z)}{P_{n+1}(0)}. \lab{bi_P} \ee In our case (taking into
account properties \re{P01} and \re{rec_coin}) it is seen that \be
\hat P_n(z) = \frac{P_{n+1}(z) - P_n(z)}{z}. \lab{bi_P_Her} \ee
Formula \re{bi_P_Her} admits another interpretation if one
introduces the Christof\/fel transform (CT) of LBP. Recall
\cite{ZheL,VZ1,VZ0} that the CT for LBP is def\/ined as \be
P_n^{(C)}(z) = \frac{P_{n+1}(z) - U_n P_n(z)}{z-\mu}, \qquad
n=0,1,\dots, \qquad U_n = \frac{P_{n+1}(\mu)}{P_n(\mu)}, \lab{CT}
\ee where $\mu$ is an arbitrary parameter such that $P_n(\mu) \ne
0$, $n=1,2,\dots$. The polynomials $P_n^{(C)}(z)$ are again monic
LBP with the transformed recurrence coef\/f\/icients \be b^{(C)}_n
= b_n \frac{b_{n+1} + U_n}{b_n + U_{n-1}}, \qquad d^{(C)}_n = d_n
\frac{d_{n+1} + U_{n+1}}{d_n + U_{n}}. \lab{bd_CT} \ee The moments
$c^{(C)}_n$ corresponding to CT are expressed as \cite{ZheL} \be
c^{(C)}_n = \frac{c_{n+1} - \mu c_n}{c_1 - \mu}. \lab{CT_moms} \ee
There is a special case of the CT when $\mu=0$. In this case, we
have \be P_n^{(C)}(z) = \frac{P_{n+1}(z) + d_n P_n(z)}{z}.
\lab{CT_0} \ee For the recurrence coef\/f\/icients in this case we
have \cite{ZheL}
\begin{gather}
\tilde b_n = b_n \frac{b_{n+1}-d_n}{b_n-d_{n-1}}, \qquad n=1,2,\dots, \nonumber \\
\tilde d_0=d_0-b_1, \qquad \tilde d_n = d_{n-1}
\frac{b_{n+1}-d_n}{b_n-d_{n-1}}, \qquad n=1,2,\dots. \lab{0CT_bd}
\end{gather}
 Note that there is some ``irregularity'' in the expression for
$\tilde d_n$ in \re{0CT_bd} for $n=0$ and $n=1,2,\dots$. This
irregularity can be avoided if one formally puts $b_0=0$. However
in our case $b_0 \ne 0$. Hence we  will indeed have such an
irregularity in the analytic dependence of the coef\/f\/icients
$\tilde d_n$ in $n$. Namely, we have the explicit expressions
\begin{gather}
\tilde b_n = -\frac{(n+1/2)^2}{(n+1)(n+2)} ,\qquad
n=1,2,\dots,\nonumber\\
 \tilde d_0 = 1/8, \qquad \tilde d_n = -\frac{n}{n+2}, \qquad
n=1,2,\dots . \lab{part_bd}
\end{gather}

The corresponding moments are transformed simply as a ``shift''
\be c_n^{(C)}=\frac{c_{n+1}}{c_1}.
 \lab{CT_c_0}
 \ee
 Comparing \re{CT_0}
with \re{bi_P_Her} we see that in our case the biorthogonal
partners $\hat P_n(z)$ coincide with the CT LBP $P_n^{(C)}$ for
$\mu=0$. Taking into account that in our case $c_1=-1$ we obtain
that the ``negative'' moments are expressed as \be c_{-n} =
-c_{n+1}, \qquad n=1,2,\dots .\lab{neg_c_Her} \ee

The LBP are connected with the two-point Pad\'e approximation
problem \cite{HR}. Given the moments $c_n$, $n=0, \pm 1, \pm 2,
\dots$, consider two formal power series \be F_{+}(z)=
\sum_{k=1}^{\infty} c_k z^{-k}, \qquad F_{-}(z)=
\sum_{k=0}^{\infty} c_{-k} z^{k}. \lab{FPM} \ee Then we have
\cite{HR,ZheL}
\begin{gather}
 \frac{P_{n-1}^{(1)}(z)}{P_n(z)}
=\frac{F_{+}(z)}{c_1} + O(z^{-n-1}), \qquad
\frac{P_{n-1}^{(1)}(z)}{P_n(z)} = -\frac{F_{-}(z)}{c_1} +
O(z^{n}). \lab{Pade1}
\end{gather}
 It is convenient to introduce the formal Laurent
series \be F(z) = \frac{F_+(z)-F_-(z)}{c_1}. \lab{F_def} \ee Then
we have \be F(z) - \frac{P_{n-1}^{(1)}(z)}{P_n(z)} = \left\{
\begin{array}{ll}
O(z^{-n-1}), \quad  & z\to \infty, \\ O(z^n), \quad & z \to
0,\end{array} \right. \lab{Pade} \ee i.e.\ LBP $P_n(z),
P_{n-1}^{(1)}(z)$ solve the problem of two-point Pad\'e
approximation (at $z=0, \infty$). Note that the Hermite formula
\re{Her_Pade} describes ``one half'' of this two-point Pad\'e
approximation problem.

So it is reasonable to refer to the polynomials $P_n(z)$ as the
Hermite elliptic Laurent biorthogonal polynomials.

\section{The weight function and biorthogonality}
 Consider the functions \be F_-(z)
=-\frac{2J(k)}{k^2 K(k)}, \qquad F_+(z) = -\frac{2
J(1/k)}{K(1/k)}, \lab{FH} \ee where $z=k^2$. It is assumed that
the function $F_-(z)$ is def\/ined near $z=0$ while the function
$F_+(z)$ is def\/ined near $z=\infty$. Note that formally
\[
F_+(z)= \frac{F_-(1/z)}{z}.
\]

From the considerations of the previous section it is easily
verif\/ied that formula \re{Pade} holds for the LBP def\/ined by
the recurrence coef\/f\/icients \re{db_Her}. It is seen also that
if one introduces the function (cf.~\cite{Hen2}) \be w(z)=
\frac{F_+(z)-F_-(z)}{2 \pi i z} \lab{weight} \ee then \be \int_C
z^{-k} w(z) dz = c_k, \qquad k=0, \pm 1, \pm 2, \dots, \lab{w_c}
\ee where the integration contour $C$ is the unit circle. Thus the
biorthogonality property \re{bi_LBP} can be presented in the form
\be \int_C P_n(z) \hat P_m(1/z) w(z)
 dz = h_{nm} \delta_{nm}. \lab{bi_w}
 \ee Now we calculate
the weight function $w(z)$ in a more explicit form. We have \be
w(z) = \frac{1}{\pi i k^2}  \left( -\frac{J(k)}{k^2 K(k)}
+\frac{J(1/k)}{K(1/k)} \right). \lab{w1} \ee Using the relation
$J(k)=K(k)-E(k)$ and the formulas for the complete elliptic
integrals with inverse modulus (as usual $K'(k)=K(k')$,
$E'(k)=E(k')$, $k'^2=1-k^2$):
\[
K(1/k)=k(K(k) + iK'(k)), \qquad E(1/k)=\frac{E(k)-iE'(k)-k'^2K(k)
+ik^2K'(k)}{k}
\] we arrive at the formula \be w(z)=
\frac{-K(k)K'(k)+E(k)K'(k)+K(k)E'(k)}{k^4 K(k) K(1/k) }.\lab{w2}
\ee The latter expression can be further simplif\/ied using the
Legendre relation \cite{WW}
\[
-K(k)K'(k)+E(k)K'(k)+K(k)E'(k)=\pi/2
\]
to give \be w(z) =  \frac{1}{2 z^{3/2}}  \frac{1}{K(k)K(1/k)}.
\lab{w3} \ee Now we introduce variable $\theta$ on the unit circle
such that $k=z^{1/2}=e^{i\theta/2}$. The biorthogonality relation
can then be written as \be \int_0^{2\pi} P_n(e^{i\theta}) \hat
P_m(e^{-i\theta}) \rho(\theta) d \theta = h_n  \delta_{nm},
\lab{bi_rho} \ee where the weight function is \be \rho(\theta) =
\frac{i}{2 e^{i\theta/2}  |K(e^{i\theta/2})|^2}. \lab{rho} \ee

\section{Generating function and explicit expression}
 From \re{AB_ser} we obtain the generating
function for the corresponding LBP \be
\Phi(x,z)=\frac{F(x;z)}{\sqrt{(1-x^2)(1-z x^2)}}
=\sum_{n=0}^{\infty} \xi_n P_n(z), \lab{gen_P} \ee where $\xi_n$
is given by \re{xi} and{\samepage
\[
F(x;z)=\int_{0}^x {\frac{d t}{\sqrt{(1-t^2)(1-z t^2)}}}
\]
is the standard (incomplete) elliptic integral of the f\/irst
kind.}

In order to f\/ind the explicit expression for the polynomials
$P_n(z)$ from \re{gen_P}, we note that if \be
\frac{1}{\sqrt{(1-x^2)(1-z x^2)}} = \sum_{n=0}^{\infty} \beta_n
x^{2n} \lab{ser_beta} \ee then, obviously,
\[
F(x;z)=\sum_{n=0}^{\infty} \beta_n x^{2n+1}/(2n+1)
\]
and hence $\Phi(x,z)=\sum\limits_{n=0}^{\infty} A_n(z) x^{2n+1}$
where \be A_n(z)=\sum_{s=0}^n {\frac{\beta_s \beta_{n-s}}{2s+1}}.
\lab{A_beta} \ee Formula \re{A_beta} gives an explicit expression
for the polynomials $A_n(z)$ and hence for the LBP $P_n(z)$ if the
coef\/f\/icients $\beta_n$ are known. But it is easy to verify
(using e.g.\ the binomial theorem) that \be \beta_n=
\frac{(1/2)_n}{n!} {_{2}F_1}(-n,1/2;1/2-n;z). \lab{beta} \ee We
thus have
\begin{gather}
A_n(z)= \sum_{m=0}^n \frac{(1/2)_n (1/2)_m (-n)_m}{(2m+1) n! m!
(1/2-n)_m} \nonumber \\
\phantom{A_n(z)=}{}\times
{_{2}F_1}(-m,1/2;1/2-m;z){_{2}F_1}(m-n,1/2;1/2+m-n;z).
\end{gather}

Another explicit expression is obtained if one notices that \be
\frac{1}{\sqrt{(1-x^2)(1-k^2x^2)}} = \sum_{n=0}^{\infty} k^n
x^{2n} Y_n\left(\frac{k+k^{-1}}{2}\right), \lab{gen_Leg} \ee where
$Y_n(z)$ are the ordinary Legendre polynomials \cite{KS}:
\[
Y_n(t) = {_{2}F_1}\left(-n,n+1;1; \frac{1-t}{2}\right).
\]
We thus have the rather simple expression{\samepage  \be A_n(z)=
k^n \sum_{s=0}^n \frac{Y_s(q)Y_{n-s}(q)}{2s+1}, \lab{P_Leg} \ee
where $q=(k+k^{-1})/2$ (recall that $k^2=z$).}

For $z=1$ we have $\beta_n =1$ for all $n$ in \re{ser_beta}. Hence
we have from \re{A_beta} \be A_n(1)=G_n, \lab{AnG} \ee where we
denote
\[
G_n = \sum_{s=0}^n \frac{1}{2s+1}
\] -- the f\/inite sum of
inverse odd numbers. $G_n$ can be obviously expressed in terms of
the Euler ``harmonic numbers'' def\/ined as
\[
H_n = \sum_{k=1}^n 1/k.
\]
We have clearly \be G_n=H_{2n+1}-H_n/2. \lab{GH} \ee

Consider the recurrence relation of the type \re{rec_AB} \be
(2n+1)\psi_{n} - 2n(1+z) \psi_{n-1} +(2n-1) z \psi_{n-2} =0.
\lab{rec_psi} \ee For $z=1$ we see that $\psi_n=G_n$ is a solution
of this equation. The second independent solution for $z=1$ is
trivial -- it is a constant: $\psi_n={\rm const}$. This means that
the general solution of the equation \re{rec_psi} for $z=1$ can be
presented in the form \be \psi_n =\alpha + \beta G_n \lab{gen_psi}
\ee with arbitrary constants $\alpha$, $\beta$. In particular, for
$B_n(1)$ we can write \be B_n(1)= G_n-1. \lab{B1} \ee

\section{Polynomials orthogonal on the unit circle}

 Consider the Christof\/fel transform
\re{CT} of our polynomials with $\mu=1$. We have \be U_n =
\frac{P_{n+1}(1)}{P_n(1)} =\frac{\xi_n  A_{n+1}(1)}{\xi_{n+1}
A_n(1)} =\frac{(n+3/2) G_{n+1}}{(n+1) G_n}. \lab{U_1} \ee For the
corresponding transformed LBP, we have the expression \be \tilde
P_n(z) =\frac{P_{n+1}(z)-U_n P_n(z)}{z-1}. \lab{tP1} \ee The
moments are calculated by \re{CT_moms}: \be \tilde c_n
=\frac{c_{n+1}-c_n}{c_1-1}=\frac{c_n-c_{n+1}}{2}. \lab{moms_1}
 \ee
Using property \re{neg_c_Her} we see that $\tilde c_{-n} = \tilde
c_n$, i.e.\ the moments $\tilde c_n$ are symmetric with respect to
ref\/lection. In turn, this is equivalent to the statement that
the corresponding polynomials $\tilde P_n(z)$ are the Szeg\"o
polynomials which are orthogonal on the unit circle
\cite{Szego,Ger}: \be \int_0^{2\pi} \tilde P_n(e^{i \theta})
\tilde P_m(e^{-i\theta}) \tilde \rho(\theta) d \theta =h_n
\delta_{nm}, \lab{PUC} \ee where \be \tilde \rho(\theta) =
\frac{\rho(\theta)}{c_1-1} (e^{i \theta} -1). \lab{c_rho} \ee In
our case we have explicitly (see \re{rho}) \be \tilde \rho(\theta)
= \frac{\sin(\theta/2)}{2|K(e^{i\theta/2})|^2}. \lab{c_rho1} \ee

It is well known  that polynomials orthogonal on the unit circle
are def\/ined by the recurrence relation \cite{Ger} \be
P_{n+1}(z)=zP_n(z)-a_n P_n^*(z), \qquad n=0,1,\dots, \lab{rec_UC}
\ee where $P_n^*(z)=z^n \bar P_n(1/z)$ (bar means complex
conjugate). In our case all moments $\tilde c_n$ are real, hence
$\tilde P_n^*(z)=z^n \tilde P_n(1/z)$. The parameters $a_n$ are
called the ref\/lection parameters. They play a~crucial role in
the theory of Szeg\"o polynomials on the unit circle. We have \be
a_n =-\tilde P_{n+1}(0). \lab{a_def} \ee It is well known
\cite{Ger} that if the ref\/lection parameters are real and
satisfy the condition $|a_n|<1$ for all $n=0,1,\dots$ then the
positive weight function $\rho(\theta)>0$ always exists. Moreover
in this case the weight function is symmetric on the unit circle:
$\rho(2\pi-\theta) = \rho(\theta)$.

In our case we have \be a_n = -\tilde P_{n+1}(0)=-U_{n+1}
P_{n+1}(0)+P_{n+2}(0) =1-U_{n+1}, \lab{aP} \ee where we used the
property $P_n(0)=1$.

It is easily seen that \be a_n = -\frac{1}{2(n+2)}
-\frac{1}{2(n+2)G(n+1)}, \qquad n=0,1,\dots. \lab{a_exp} \ee In
\re{a_exp} both terms are negative and less then 1/2 in absolute
value. Hence $-1<a_n<0$ for all $n=0,1,\dots$ and the function
$\tilde \rho(\theta)$ is positive as is seen from \re{c_rho1}.

We thus have a (presumably) new example of the Szeg\"o polynomials
orthogonal on the unit circle for which both weight function and
recurrence coef\/f\/icients are known explicitly.

Following \cite{DG} and \cite{ZheC}, to any polynomials $\tilde
P_n(z)$ orthogonal on the unit circle with the property $-1 < a_n
<1$ one can associate symmetric monic polynomials $S_n(x)=x^n+
O(x^{n-1})$ orthogonal on an interval of the real axis. Explicitly
\be S_n(x) = \frac{z^{-n/2}(\tilde P_n(z) + \tilde
P^*_n(z))}{1-a_{n-1}}, \lab{DGS} \ee where $x=z^{1/2} + z^{-1/2}$
(it is assumed that one chooses one branch of the function
$z^{1/2}$ such that for $z=r  e^{i \theta}$ we have
$z^{1/2}=r^{1/2} e^{i\theta/2}$, $ -\pi < \theta< \pi$). The
polynomials $S_n(x)$ satisfy the three-term recurrence relation
\be S_{n+1} + u_n S_{n-1}(x) = xS_n(x), \lab{rec_S} \ee where the
recurrence coef\/f\/icients are \be u_n = (1+a_{n-1})(1-a_{n-2}),
\qquad n=1,2,\dots. \lab{uDG} \ee In \re{uDG} it is assumed that
$a_{-1}=-1$ (this is a standard convention in the theory of
polynomials orthogonal on the unit circle), so $u_1=2(1+a_0)$.

If the polynomials $\tilde P_n(z)$ are orthogonal on the unit
circle \be \int_{0}^{2\pi} \tilde P_n(e^{i\theta}) \tilde
P^*_m(e^{-i \theta}) \rho(\theta) d\theta = 0, \qquad m \ne n
\lab{ort_c} \ee with the weight function $\rho(\theta)$ then
polynomials $S_n(x)$ are orthogonal on the symmetric interval
$[-2,2]$ \be \int_{-2}^2 S_n(x) S_m(x) w(x) dx = h_n \delta_{nm}
\lab{ort_S} \ee with the weight function \cite{DG,ZheC} \be
w(x)=\frac{\rho(\theta)}{\sin(\theta/2)}, \lab{w_rho} \ee where
$x=2\cos(\theta/2)$.

In our case it is elementary verif\/ied that \be u_n =
\frac{(n+1/2)^2}{n(n+1)}, \qquad n=1,2,\dots. \lab{u_Leg} \ee
These recurrence coef\/f\/icients correspond to well known
associated Legendre polynomials studied, e.g.\ in \cite{BD}.
Recall that generic associated Legendre polynomials are symmetric
OP satisfying the recurrence relation \re{rec_S} with the
recurrence coef\/f\/icients \be u_n =
\frac{(n+\nu)^2}{(n+\nu)^2-1/4} \lab{Leg_nu} \ee with arbitrary
nonnegative parameter $\nu$. The ordinary Legendre polynomials
correspond to $\nu=0$. In our case we have $\nu=1/2$.

Consider the weight function for these polynomials. From
\re{w_rho} and \re{c_rho1} we derive \be w(x) =
\frac{1}{2|K(e^{-i\theta/2})|^2}. \lab{wL} \ee We can simplify
this formula if one exploits the relations for $K(z)$ where
$|z|=1$. Namely, one has \cite{MagOb}
\[
K(e^{\pm i \phi}) = \frac{1}{2} e^{\mp i \phi/2} (K(\cos(\phi/2))
\pm i K(\sin(\phi/2))), \qquad -\pi < \phi \le \pi,
\]
whence we have \be w(x)= \frac{2}{K^2(\cos(\theta/4)) +
K^2(\sin(\theta/4))}. \lab{wL2} \ee Taking into account that $x=2
\cos(\theta/2)$ we f\/inally arrive at the formula \be w(x)=
\frac{2}{K^2(\sqrt{1/2+x/4}) + K^2(\sqrt{1/2-x/4})}, \qquad -2 \le
x \le 2. \lab{wL3} \ee Note that the function $w(x)$ is even
$w(-x)=w(x)$ as should be for symmetric polynomials. It has the
only maximum at $x=0$: $w(0)=1/K^2(\sqrt{1/2})=\frac{16
\pi}{\Gamma^4(1/4)}$. Near the endpoints of the interval $[-2,2]$
the weight function $w(x)$ tends to zero rapidly.

The weight function for generic associated Legendre polynomials
(with arbitrary $\nu$) was found in \cite{BD}. Our formula
\re{wL3} can be obtained from the results of \cite{BD} by putting
$\nu=1/2$. Note, that Pollaczek studied \cite{Poll} more general
orthogonal polynomials containing 4 parameters. The associated
Legendre polynomials (as well as the associated ultraspherical
polynomials) are contained in the Pollaczek polynomial family as a
special case.

Similar polynomials were studied also in \cite{Rees} where the
author in fact rediscovered the Hermite approach (as well
Hermite's continued fraction \re{Her_cf}) to elliptic integrals.
He considered integrals $J_n(k)$ as moments for some ``elliptic''
orthogonal polynomials having $w(x)=1/\sqrt{(1-x^2)(1-k^2x^2)}$ as
an orthogonality weight on the interval $[-1,1]$. As a by-product
the author of \cite{Rees} introduced polynomials which are similar
to the Hermite polynomials $A_n(z)$, $B_n(z)$. He then related
them with the associated Legendre polynomials in a similar way. It
is interesting to note that Hermite himself already introduced
such polynomials in~\cite{Her1,Her2}. Hermite also established
dif\/ferential properties of these polynomials anticipating
results of Rees~\cite{Rees}.

Consider relations between LBP of special type and orthogonal
polynomials on an interval in details.

We have $\tilde P_n(z) = (P_{n+1}(z)-U_nP_n(z))/(z-1)$.
Substituting this formula to \re{DGS} and using the inversion
property $z^nP_n(1/z)=P_n(z)$ of our LBP $P_n(z)$ we immediately
obtain a very simple relation \be S_n(x) = z^{-n/2} P_n(z).
\lab{SP} \ee We can see this also using recurrence relation for
the LBP \be P_{n+1}(z) + d_n P_n(z) = z(P_n(z) + b_n P_{n-1}(z)).
\lab{recLBP} \ee Assume (as in our case) that $d_n=-1$,
$n=0,1,\dots$. Then $P_n(0)=1$ for all $n=0,1,2,\dots$ and from
\re{rec_bd} we obtain that the reciprocal polynomials coincide
with the initial ones: \be z^n P_n(1/z) = P_n(z). \lab{rec_pr} \ee
Conversely, assume that some LBP are reciprocal \re{rec_pr}. Then
from \re{rec_bd} it follows that either $d_n=1,\; n=0,1,\dots$ or
$d_n=-1, \; n=0,1,\dots$. But we have obviously $P_n(0)=1$ which
leads to the only possibility $d_n=-1, \; n=0,1,2,\dots$. Thus
condition $d_n=-1$ for all $n$ is necessary and suf\/f\/icient for
LBP to be reciprocal invariant \re{rec_pr}. In this case
polynomials \be S_n(x) = z^{-n/2} P_n(z) \lab{SP_rec} \ee are
obviously monic polynomials in $x=z^{1/2} + z^{-1/2}$. From
recurrence relation \re{recLBP} with $d_n=-1$ we obtain recurrence
relation for polynomials $S_n(x)$ \be S_{n+1}(x) + u_n S_{n-1}(x)
= x S_n(x), \lab{rec_S1} \ee where $u_n=-b_n$. Thus we arrived at
the same symmetric polynomials on the interval as in the case of
polynomials orthogonal on the unit circle. Note that for the LBP
$P_n(z)$ with the property \re{rec_pr} we have $c_1=d_0=-1$ and
$c_n^*=c_n$ thus from \re{mom_rec} we obtain \be c_{1-n}=-c_n
\lab{sym_mom} \ee for all $n=0,\pm 1, \pm 2,\dots$. Perform now
the Christof\/fel transform with $\mu =1$ \be \tilde P_n(z) =
\frac{P_{n+1}(z) - U_n P_n(z)}{z-1}. \lab{CT1} \ee From
\re{CT_moms} and \re{sym_mom} we obtain that transformed moments
are symmetric $\tilde c_{-n}=c_n$. This means that polynomials
$\tilde P_n(z)$ will satisfy the recurrence relation \re{rec_UC}
with ref\/lection parameters $a_n$ given by $a_n=1-U_{n+1} =
1-P_{n+2}(1)/P_{n+1}(1).$ Formula~\re{DGS} now is equivalent to
formula~\re{SP_rec}. Thus starting from arbitrary LBP with the
property $d_n=-1$ we can arrive at the same symmetric OP on the
interval $S_n(x)$ as for the case of polynomials orthogonal on the
unit circle. From another point of view these relations are
discussed also in \cite{Ranga}.

\section{Geronimus transform. Laurent biorthogonal polynomials\\ with a
concentrated mass added to the measure}

In the previous section we showed that the Christof\/fel
transformation of Hermite's elliptic LBP gives polynomials
orthogonal on the unit circle with explicit ref\/lection
coef\/f\/icients \re{a_exp} and the weight function given by
\re{c_rho1}. In this section we consider another spectral
transformation of Hermite's elliptic LBP which is called the
Geronimus transform (GT). This transform was introduced in
\cite{ZheL} and is similar to well known Geronimus transform for
the ordinary orthogonal polynomials \cite{SVZ,ZheR}. Recall thatGT
for LBP is def\/ined as \cite{ZheL} \be \tilde P_n(z) = V_n P_n(z)
+z(1 - V_n )P_{n-1}(z), \lab{GT} \ee where \be V_0=1, \qquad V_n =
\frac{\mu}{\mu - \phi_n/\phi_{n-1}}, \qquad n=1,2,\dots. \lab{VG}
\ee In \re{VG} $\mu$ is an arbitrary parameter and $\phi_n$ is an
arbitrary solution of the recurrence relation \be \phi_{n+1} + d_n
\phi_n = \mu (\phi_n + b_n \phi_{n-1}). \lab{phi_rec} \ee Note
that \re{phi_rec} is the same recurrence relation that
\re{rec_LBP} for LBP. Hence its the general solution (up to a
common factor) can be presented in the form \be \phi_n = P_n(\mu)
+ \chi P_{n-1}^{(1)}(\mu). \lab{phi_sol} \ee It is easy to verify
that the polynomials $\tilde P_n(z)$ are again LBP satisfying
recurrence relation \be \tilde P_{n+1}(z) + \tilde d_n \tilde
P_{n}(z) = z(\tilde P_{n}(z) + \tilde b_n \tilde P_{n-1}(z)),
\lab{rec_G} \ee where the recurrence coef\/f\/icients are \be
\tilde b_1 = \chi (V_1-1), \qquad \tilde b_n = b_{n-1}
\frac{1-V_n}{1-V_{n-1}}, \qquad \tilde d_n = d_n
\frac{V_{n+1}}{V_n}. \lab{bd_G} \ee It can be shown that CT and GT
are reciprocal to one another~\cite{ZheL}. This observation allows
one to obtain explicit biorthogonality relation for polynomials
$\tilde P_n(z)$ starting from that for polynomials $P_n(z)$.
Namely, the pair of the Stieltjes functions $F_+(z)$, $F_-(z)$ is
transformed as~\cite{ZheL} \be \tilde F_+(z)= \frac{\nu F_+(z) +
\mu + \nu}{z- \mu}, \qquad F_-(z)= \frac{\nu F_-(z) -\mu - \nu}{z-
\mu}, \lab{FF_G} \ee where
\[
\nu = \frac{\mu \chi}{c_1-\chi} = \frac{\mu \chi}{d_0-\chi}.
\]
Assume that $|\mu| \le 1$. Choose the contour $C$ as the unit
circle (if $|\mu|=1$ we choose $C$ as the unit circle with a small
deformation near $z=\mu$ in order to include the point $z=\mu$
inside the contour $C$). Then the weight function for new
polynomials $\tilde P_n(z)$ is def\/ined as in \re{weight}, i.e.
\be \tilde w(z) = \frac{\tilde F_+(z) - \tilde F_-(z)}{2\pi i z} =
\frac{w(z)}{z- \mu} +   \frac{2(\mu +\nu)}{2 \pi i z (z-\mu)}.
\lab{w_G} \ee The second term in \re{w_G} will give a concentrated
mass at the point $z=\mu$ added to a ``regular'' part presented by
the f\/irst term in \re{w_G}.

Assume that $\mu=1$, as in our case. Then we have biorthogonality
relation \be \int_0^{2\pi} \tilde \rho(\theta) \tilde
P_n(e^{i\theta}) {\hat {\tilde  P}}_n(e^{-i\theta}) d \theta + M
\: \tilde P_n(1) {\hat {\tilde  P}}_n(1) = 0, \qquad n \ne m,
\lab{bi_G} \ee where \be \tilde \rho(\theta) =
\frac{\rho(\theta)}{e^{i\theta}-1} \lab{reg_rho} \ee is the
``regular'' of the weight function on the unit circle and the last
term \re{reg_rho} describes the concentrated mass \be M= \frac{\nu
+1}{2 \pi} \lab{M_G} \ee inserted at the point $z=1$.

In our case we have $d_n =-1$, $n=0,1,\dots$ hence
$\nu=-\chi/(\chi+1)$ and
\[
M=\frac{1}{2\pi (1+\chi)}.
\] For the
recurrence coef\/f\/icients we have explicit formulas \re{bd_G}
where we need f\/irst to calculate the coef\/f\/icients $V_n$. We
see that
\[
\phi_n = P_n(1) + \chi P_{n-1}^{(1)}(1).
\] But the values $P_n(1)$
and $P_{n-1}^{(1)}(1)$ were already calculated (see \re{AnG},
\re{B1}). We thus have \be \phi_n = \frac{(3/2)_n }{n!}(G_n + 2
\chi (G_n-1)). \lab{phi_H} \ee Now all coef\/f\/icients $V_n$,
$\tilde b_n$, $\tilde d_n$ are calculated explicitly.

We thus constructed a nontrivial example of the Laurent
biorthogonal polynomials with explicit both recurrence
coef\/f\/icients and the measure. The weight function for these
polynomials has a concentrated mass on the unit circle. Note that
in contrast to the Christof\/fel transformed Hermite's polynomials
the polynomials $\tilde P_n(z)$ constructed in this section are
not polynomials of the Szeg\"o type. This means, in particular,
that the biorthogonal partners ${\hat {\tilde P}}_n(z)$ do not
coincide with $\tilde P_n(z)$.

\section{Associated families of the Laurent biorthogonal
polynomials}

 Return to the sequence
$J_n(x,k)$ of incomplete elliptic integrals def\/ined by \re{J_n}.
They satisfy three-term recurrence relation \re{J_n_rec}.
Repeating previous considerations we can express $J_n(x;k)$ for
any $n=1,2,\dots$ in terms of $J_j(x;k)$ and $J_{j+1}(x;k)$ for
some f\/ixed nonnegative integer $j$ (in~\re{JnQ} the case $j=0$
is chosen): \be J_{n+1}(x;k)= Q(x;n,j)  \sqrt{(1-x^2)(1-k^2x^2)} +
A_{n-j}^{(j)} J_{j+1}(x;k) - B_{n-j}^{(j)} J_j(x;k)\lab{Jnj} \ee
with some coef\/f\/icients $A_n^{(j)}$, $B_n^{(j)}$.

We f\/irst note that there is an explicit expression of complete
elliptic integrals $J_n(k)$ in terms of the Gauss hypergeometric
function \be J_n(k)= \frac{k^{2n} \pi (1/2)_n }{2 n!} \,
{_{2}F_1}(1/2,1/2+n;1+n;k^2) \lab{Jn_hyp} \ee (relation
\re{Jn_hyp} can be easily verif\/ied by direct integration).

Repeating similar considerations that were already exploited in
the two f\/irst sections we can show that $A_n^{(j)}(z)$ and
$B_n^{(j)}$ are determined by the recurrence relations
\begin{gather}
(n+3/2+j) A_{n+1}^{(j)} =(n+j+1)(z+1)A_n^{(j)} -
z(n+j+1/2)A_{n-1}^{(j)}, \nonumber \\
(n+3/2+j) B_{n+1}^{(j)} =(n+j+1)(z+1)B_n^{(j)} -
z(n+j+1/2)B_{n-1}^{(j)} \lab{jAB}
\end{gather}
with the same initial conditions as \re{ini_AB}. Hence,
$A_n^{(j)}(z)$ polynomials of $n$-th degree in $z=k^2$ and
$B_n^{(j)}(z)$ is a polynomial of degree $n-1$ multiplied by $z$.

Introduce also the function \be F(z;j) =
\frac{J_{j+1}(k)}{J_j(k)}, \lab{Fj} \ee where as usual we put
$k^2=z$.

Then we will have the property \be F(z;j) -
\frac{B_n^{(j)}(z)}{A_n^{(j)}(z)} =O(z^{n+1}). \lab{pade_j} \ee

In the same way it is possible to derive the generating functions
for polynomials $A_n^{(j)}(z)$, $B_n^{(j)}(z)$:
\begin{gather}
(2j+1)\frac{k^{-2j} J_j(x;k)}{\sqrt{(1-x^2)(1-k^2x^2)}} =
\sum_{n=j}^{\infty} A_{n-j}^{(j)}(z) x^{2n+1}, \nonumber\\
(2j+1)\frac{k^{-2j}J_{j+1}(x;k)}{\sqrt{(1-x^2)(1-k^2x^2)}} =
\sum_{n=j}^{\infty} B_{n-j}^{(j)}(z) x^{2n+1}. \lab{assoc_gen}
\end{gather}
Using \re{gen_Leg} we arrive at the explicit representation for
polynomials $A_n^{(j)}(z)$ in terms of the Legendre polynomials
$Y_n(q)$: \be A_n^{(j)}(z) = (2j+1) k^n  \sum_{s=0}^n
\frac{Y_s(q)Y_{n-s}(q)}{2s+2j+1}, \lab{Aj_exp} \ee where $z=k^2$,
$q=(k+1/k)/2$.

Now we introduce the monic LBP $P_n^{(j)}(z) = A_n^{(j)}(z)/\xi_n$
where $\xi_n = (j+1)_n/(j+3/2)_n$. They satisfy the recurrence
relation of type \re{rec_LBP} with \be d_n =-1, \qquad b_n
=-\frac{(n+j+1/2)^2}{(n+j)(n+j+1)}. \lab{ass_recP} \ee Thus we
have $j$-associated polynomials with respect to the Hermite
elliptic LBP, i.e.\ we should replace $n \to n+j$ in formulas for
recurrence coef\/f\/icients.

Due to the property $d_n =-1$ the associated polynomials
$P_n^{(j)}$ again possess the invariance property
$z^nP_n^{(j)}(1/z)=P_n^{(j)}(z)$. We thus can construct
polynomials orthogonal on the unit circle using formula \re{CT1}.
In order to get explicit expression for ref\/lection parameters
$a_n^{(j)}$ we need the value $P_n^{(j)}(1)$. But this can easily
be obtained from \re{Aj_exp}: \be A_n^{(j)}(1)=  G_n(j) \lab{Aj1}
\ee where we introduced the function
\[
G_n(j)=(2j+1)  \sum_{s=0}^n \frac{1}{2s+2j+1}.
\]
(Up to a common factor $G_n(j)$ is a sum of $n+1$ succeeding
inverse odd numbers starting from $1/(2j+1)$; for $j=0$ it
coincides with $G_n$.)
 Thus
\[
U_n=\frac{P_{n+1}^{(j)}(1)}{P_{n}^{(j)}(1)} =
\frac{n+3/2+j}{n+j+1} \: \frac{G_n(j)}{G_{n+1}(j)}
\]
and
\[
a_n^{(j)}=1-U_{n+1}.
\]
Corresponding symmetric polynomial $S_n^{(j)}(x)$ on the interval
satisfy recurrence relation \re{rec_S1} with \be
u_n=\frac{(n+j+1/2)^2}{(n+j)(n+j+1)} \lab{assL} \ee so they
coincide with the associated Legendre polynomials considered
in~\cite{BD}. Indeed, from formula \re{Leg_nu} we see that the
``shift'' parameter $\nu = j+1/2$ where $j=0,1,2,\dots$.

It is interesting to note that recurrence relations for
polynomials $P_n^{(j)}$ can be presented in such a form  that the
all coef\/f\/icients are linear in $n$ -- see e.g. \re{jAB}. In
\cite{GVZ} we considered a family of such LBP connected with
so-called generalized eigenvalue problem on $su(1,1)$ Lie algebra.
In our case, however, corresponding representations of $su(1,1)$
will not be unitary, in contrast to~\cite{GVZ}. This leads to an
interesting open problem how to describe associated classical LBP
in terms of non-unitary representations of $su(1,1)$ algebra. Note
that generic orthogonal polynomials with linear recurrence
coef\/f\/icients in $n$ were studied in details by Pollaczek
\cite{Poll} who derived explicit expression for them and found the
weight function as well.

Note f\/inally that LBP with recurrence coef\/f\/icients
\re{ass_recP} belong to a family of so-called ``associated Jacobi
Laurent polynomials'' introduced and studied by Hendriksen
\cite{Hen1, Hen2}. Nevertheless, in our case the parameters of the
associated LBP belong to the exceptional class which was not
considered in \cite{Hen1, Hen2}. This means that in some formulas
in \cite{Hen1, Hen2} the bottom parameter in the Gauss
hypergeometric  function $_{2}F_1(z)$ takes negative integer
values. In this case formulas obtained by Hendriksen should be
rederived in a dif\/ferent form. In particular, a simple explicit
expression of the power coef\/f\/icients (in terms of the
hypergeometric function $_{4}F_3(1)$) for associated LBP obtained
in \cite{Hen1} seems not to be valid in our case. Instead, we
obtained explicit expressions like \re{Aj_exp}.

\section[Laurent biorthogonal polynomials connected with the
Stieltjes-Carlitz elliptic polynomials]{Laurent biorthogonal
polynomials connected\\ with the Stieltjes--Carlitz elliptic
polynomials}

Return to formula \re{sn_rec}. If one denotes \be
r_n(p)=\int_0^{\infty} k^{2n} sn^{2n}(t) e^{-pt} dt \lab{Lap_sn}
\ee then we obtain the recurrence relation \be p^2 r_n = 2n(2n+1)
r_{n+1} -4(1+k^2) n^2 r_n + 2n(2n-1) r_{n-1}. \lab{rec_r} \ee
Again it is seen that for every $n>0$ one can present \be r_{n+1}
= A_n r_1 - B_n r_0, \lab{r_nAB} \ee where obviously
$r_0=\int_0^{\infty} e^{-pt} dt = p^{-1}$ and  $A_n$, $B_n$
satisfy the same recurrence relations as $r_{n+1}$ e.g. \be p^2
A_n = 2(n+1)(2n+1) k^2 A_{n-1} -4(1+k^2) (n+1)^2 A_n +
2(n+1)(2n+3) A_{n+1} \lab{rec_AL} \ee with initial conditions
$A_0=1$, $A_{-1}=0$, $B_0=0$, $B_{-1}=-1$. Now it is seen that
$A_n$ are polynomials of degree $n$ {\it in both} variables $p^2$
and $k^2$. If $p=0$ then $A_n$ become polynomials in $k^2$
introduced by Hermite.

Coef\/f\/icients $A_n$ considered as polynomials in $p^2$ become
{\it orthogonal polynomials} because they satisfy three-term
recurrence relation typical for orthogonal polynomials. Orthogonal
polynomials of such type (and several related ones) where
introduced and studied by Carlitz. He exploited some explicit
continued fractions found by Stieltjes. These continued fractions
are connected with elliptic functions (for details see,
e.g.~\cite{IM}). Today these orthogonal polynomials are known by
Stieltjes--Carlitz elliptic polynomials \cite{Bril,Chi,IM}. Note
that the Stieltjes continued fraction is obtained from \re{rec_r}
by the same way as Hermite obtained his continued fraction
\re{Her_cf} for a ratio of two elliptic integrals. Consider now
polynomials $A_n$ as LBP with respect to the argument $z=k^2$.
Passing from $A_n$ to monic polynomials $P_n(z)$ (by the same way
as for the case $p=0$) we arrive at the recurrence relation
\re{rec_LBP} with \be d_n =-1 -\frac{p^2}{4(n+1)^2}, \qquad
b_n=-\frac{(n+1/2)^2}{n(n+1)}. \lab{bd_SC} \ee We see that the
recurrence coef\/f\/icient~$b_n$ is the same as for the Hermite
LBP, but the coef\/f\/icient~$d_n$ now depends on~$n$. This means
that polynomials $P_n(z)$ do not possess symmetric property
like~\re{rec_coin}.

In contrast to the case $p=0$ the polynomials $P_n(z)$ have more
complicated properties. For example, the reciprocal polynomials
$P_n^*(z)$ def\/ined by \re{rec_P_def} do not belong to the same
class, their recurrence coef\/f\/icients appear to be \be d^*_n=
-\frac{1}{1+\frac{p^2}{4(n+1)^2}}, \qquad b^*_n =
-\frac{(n+1)(n+2)^2(n+1/2)^2}{((n+2)^2+p^2/4)((n+1)^2+p^2/4)}.
\lab{bd_r_SC} \ee The biorthogonal partners $\hat P_n(z)$
def\/ined by \re{bi_P} have the recurrence coef\/f\/icients
\begin{gather} \hat d_n =
-\frac{n(n+1)(p^2(n+1) -n-2)}{(n(p^2-1) -1)((n+2)^2+p^2/4)}, \nonumber\\
\hat b_n= - \frac{(n+1/2)^2 (n+1)^2 (p^2(n+1) -n-2)}{(n(p^2-1)
-1)((n+1)^2 +p^2/4)((n+2)^2 +p^2/4)}. \lab{bd_bi_SC}
\end{gather} Hence
biorthogonal partners $\hat P_n(z)$ also do not belong to the same
class that initial LBP $P_n(z)$.

For the ordinary Stieltjes--Carlitz elliptic orthogonal
polynomials the orthogonality measure can be found explicitly: it
is purely discrete one (see \cite{Chi,IM}). The problem of
f\/inding the orthogonality measure for the corresponding LBP
seems to be much more complicated.

\subsection*{Acknowledgments}

The authors thank to the referees for their remarks leading to
improvement of the text. A.Zh. thanks Centre de Recherches
Math\'ematiques of the Universit\'e de Montr\'eal for hospitality.

\newpage

\pdfbookmark[1]{References}{ref}
\begin{thebibliography}{99} 

\footnotesize\itemsep=0pt

\bi{BD} Barrucand P., Dickinson D., On the associated Legendre
polynomials, in Orthogonal Expansions and Their Continued
Analogues, Editor D.T.~Haimo, Southern Illinois Press, 1968,
43--50.

\bi{Ranga} Bracciali C.F., da Silva  A.P., Sri Ranga A., Szeg\"o
polynomials: some relations to $L$-ortho\-gonal and ortho\-gonal
polynomials, {\it J. Comput. Appl. Math.} {\bf 153} (2003), 79--88.

\bi{Chi} Chihara T., An introduction to orthogonal polynomials,
Gordon and Breach, 1978.

\bi{DG} Delsarte P.,  Genin Y., The split Levinson algorithm, {\it
IEEE Trans. Acoust. Speech Signal Process} {\bf 34} (1986),
 470--478.

\bi{Ger} Geronimus Ya.L., Polynomials orthogonal on a circle and
their applications, in {\it Amer. Math. Soc. Transl. Ser. 1},
Vol.~3, Amer. Math. Soc., Providence,   1962, 1--78.

\bi{GVZ} Gr\"unbaum F.A., Vinet L., Zhedanov A., Linear operator
pencils on Lie algebras and Laurent biorthogonal polynomials, {\it
J. Phys. A: Math. Gen.} {\bf 37} (2004), 7711--7725.

\bi{HR} Hendriksen E., van Rossum H., Orthogonal Laurent
polynomials, {\it Indag. Math. (Ser. A)} {\bf 48} (1986), 17--36.

\bi{Hen1} Hendriksen E., Associated Jacobi--Laurent polynomials,
{\it J. Comput. Appl. Math.} {\bf 32} (1990), 125--141.

\bi{Hen2} Hendriksen E., A weight function for the associated
Jacobi--Laurent polynomials, {\it J. Comput. Appl. Math.} {\bf 33} (1990), 171--180.

\bi{Her1} Hermite C., Sur la d\'eveloppement en s\'erie des
integrales elliptiques de premiere et de seconde espece, {\it
Annali
di Matematica} {\bf II}  (1868), 2 ser., 97--97.

\bi{Her2}
Hermite C., Oeuvres, Tome II, Paris, 1908, 486--488.

\bi{Her} Hermite C.,  Cours d'analyse de la Facult\'e des
Sciences, Editor Andoyer,  Hermann, Paris, 1882 (Lithographed
notes).

\bi{IM} Ismail M.E.H., Masson D., Some continued fractions related
to elliptic functions, {\it Contemp. Math.} {\bf 236} (1999),
149--166.

\bi{JT} Jones W.B., Thron W.J., Survey of continued fraction
methods of solving moment problems, in Analytic Theory of
Continued Fractions, {\it Lecture Notes in Math.}, Vol.~932,
 Springer, Berlin~-- Heidelberg~-- New York,  1981.

\bi{KS} Koekoek R., Swarttouw R.F.,  The Askey scheme of
hypergeometric orthogonal polynomials and its $q$-analogue, Report
94-05, Delft University of Technology, Faculty of Technical
Mathematics and Informatics,
 1994.

\bi{Bril} Lomont J.S., Brillhart J., Elliptic polynomials,
Chapman \& Hall/CRC,  Boca Raton, FL, 2001.

\bi{MagOb} Magnus W., Oberhettinger F., Formeln und S\"atze fur
die speziellen Functionen der mathematischen Physik, Springer, Berlin,
 1948.

\bi{Poll} Pollaczek F., Sur une famille de polyn\^omes orthogonaux
\`a quatre paramitres, {\it C. R. Acad. Sci. Paris} {\bf 230} (1950),
2254--2256.

\bi{Rees} Rees C.J., Elliptic orthogonal polynomials, {\it Duke
Math. J.} {\bf 12}  (1945), 173--187.

\bi{SVZ} Spiridonov V., Vinet L., Zhedanov A., Spectral
transformations, self-similar reductions and orthogonal
polynomials, {\it J. Phys. A: Math. Gen.} {\bf 30} (1997), 7621--7637.

\bi{Szego}Szeg\"o G., Orthogonal polynomials, AMS, 1959.

\bi{VZ0} Vinet L., Zhedanov A., An integrable chain and
bi-orthogonal polynomials, {\it Lett. Math. Phys.} {\bf 46} (1998),
233--245.

\bi{VZ1} Vinet L., Zhedanov A., Spectral transformations of the
Laurent biorthogonal polynomials. I. $q$-Appel polynomials, {\it
J. Comput. Appl. Math.} {\bf 131} (2001), 253--266.

\bi{WW} Whittacker E.T., Watson G.N., A course of modern analysis,
4th ed., Cambridge University Press, 1927.

\bi{ZheR} Zhedanov A., Rational spectral transformations and
orthogonal polynomials, {\it J. Comput. Appl. Math.} {\bf 85} (1997),
67--86.

\bi{ZheC} Zhedanov A., On some classes of polynomials orthogonal
on arcs of the unit circle connected with symmetric orthogonal
polynomials on an interval, {\it J. Approx. Theory} {\bf 94} (1998),
73--106.

\bi{ZheL} Zhedanov A., The ``classical'' Laurent biorthogonal
polynomials, {\it  J. Comput. Appl. Math.} {\bf 98} (1998), 121--147.
\eb\LastPageEnding

\end{document}